\documentclass{article}[12pt]

\usepackage{amssymb, amsmath, amsthm,commath,comment,stackrel,geometry,enumerate,mathtools,bbold,xcolor,thm-restate}

\author{Jeroen Winkel}
\date{}
\title{Cycles in graphs with geometric property (T)}

\theoremstyle{definition}
\newtheorem{example}{Example}[section]
\theoremstyle{plain}

\newtheorem{lemma}[example]{Lemma}
\newtheorem{proposition}[example]{Proposition}

\newtheorem{corollary}[example]{Corollary}

\theoremstyle{remark}
\newtheorem{remark}[example]{Remark}
\theoremstyle{definition}
\newtheorem{definition}[example]{Definition}

\renewcommand{\phi}{\varphi}
\renewcommand{\H}{\mathcal H}
\renewcommand{\epsilon}{\varepsilon}
\DeclareMathOperator{\cs}{cs}

\DeclareMathOperator{\tensor}{\otimes}

\begin{document}
\maketitle

\begin{abstract}
    We show that a sequence of graphs with uniformly bounded vertex degrees, number of vertices going to infinity, and with geometric property (T) has many small cycles.
    We also show that when a small part of such a sequence of graphs with geometric property (T) is changed, it still has geometric property (T), provided that it is still an expander.
    We use this to give an example of a sequence of graphs with geometric property (T) that has large cycle-free balls.
\end{abstract}

\section{Introduction}
In this paper, we are interested in sequences $(X_n)$ of finite connected graphs.
All graphs we consider are simple and undirected.
The vertex sets will be denoted by $V(X_n)$ and the edge sets by $E(X_n)$.
We will always assume that $\lim_{n\to\infty}|V(X_n)|=\infty$ and that the sequence has \emph{uniformly bounded degree}, i.e. there is a constant $d$ such that for all $n$, all vertices of $X_n$ have degree at most $d$.

Let $X$ be a finite graph.
For vertices $x,y$ we write $x\sim y$ if $(x,y)$ is an edge, and we write $\deg(x)$ for the degree of $x$.
The \emph{Laplacian} of $X$ is the matrix $\Delta_X$, whose rows and columns are indexed by the vertices of $X$, defined by
\[(\Delta_X)_{xy} = \begin{cases}\deg(x)&\text{ if }x=y\\-1&\text{ if }x\sim y\\0&\text{ else.}\end{cases}\]
It is a symmetric positive semi-definite matrix, and it has 0 has an eigenvalue.
This is a simple eigenvalue if and only if $X$ is connected.

Consider a sequence $X=(X_n)$ of finite connected graphs with uniformly bounded degree and number of vertices going to infinity.
The sequence is an \emph{expander sequence} if there is a constant $h>0$ such that for all $n$, all positive eigenvalues of $\Delta_{X_n}$ are at least $h$.
Equivalently, we can directly look at the operator
\[\Delta_X = \bigoplus_n\Delta_{X_n} \subseteq B(L^2X).\]
Here and in the following $L^2X$ denotes the Hilbert space generated by the vertices of $X_1,X_2,\ldots$ as basis vectors.
Then the sequence is an expander sequence if and only if there is $h>0$ such that $\sigma(\Delta_X)\subseteq \{0\}\cup [h,\infty)$.

There is an equivalent combinatorial characterisation of expanders.
For a subset of the vertices $A\subseteq V(X_n)$, let $\delta A$ denote the set of edges with exactly one vertex in $A$.
Then the sequence $(X_n)$ is an expander if and only if it has uniformly bounded degree and there is $c>0$ such that for all $A\subset V(X_n)$ with $|A|\leq \frac12|V(X)|$, we have $|\delta A| \geq c|A|$.
We refer to e.g. \cite{Lub94} for more details on expanders.

In \cite{WilYu12}, Willett and Yu introduced Geometric property (T).
It was studied in more depth by the same authors in \cite{WilYu14}.
It is a stronger property than being an expander, based on spectral gap of $\Delta_X$ in a larger algebra.

Let $(X_n)$ be a sequence of graphs and let $T$ be a bounded operator in $\prod_nB(L^2X_n) \subseteq B(L^2X)$.
For $x,y\in V(X_n)$, we denote by $d(x,y)$ the shortest-path distance from $x$ to $y$.
We can view $T$ as a matrix whose rows and columns are indexed by the vertices of $X$, and $T_{xy}=0$ if $x\in X_n,y\in X_m,m\neq n$.
The \emph{propagation} of $T$ is $\sup\{d(x,y) \mid T_{xy}\neq 0\}$, which may be infinite.
The sum and product of operators with finite propagation have finite propagation again.
The \emph{algebraic uniform Roe algebra}, as introduced by Roe in \cite{Roe93}, is defined as
\[\mathbb C_{\cs}[X]  = \left\{T\in \prod_nB(L^2X_n) \mid T\text{ has finite propagation}\right\}.\]

The algebraic uniform Roe algebra is a unital pre-$C^*$-algebra.
A representation of $\mathbb C_{\cs}[X]$ is given by a Hilbert space $\H$ and a unital $*$-homomorphism $\pi\colon \mathbb C_{\cs}[X]\to B(\H)$.

The completion of $\mathbb C_{\cs}[X]$ in the norm inherited by $B(L^2X)$ is the \emph{reduced uniform Roe algebra} $C^*_{\text{red}}(X)$.
In this paper we are more interested in the completion with respect to a larger norm, namely the maximal norm $\norm{\cdot}_{\max}$ given by
\[\norm T_{\max} = \sup_{(\pi,\H)} \norm{\pi(T)},\]
where the supremum is taken over all representations $(\pi,\H)$.
If $X$ has uniformly bounded degree, then the maximal norm is necessarily finite (see for example \cite[Section 3]{GonWanYu08}).
It is then easy to check that it is indeed a norm.
Its completion is a $C^*$-algebra, the \emph{maximal uniform Roe algebra} $C^*_{\max}(X)$ (see for example \cite{WilYu12I}).

For $T\in \mathbb C_{\cs}[X]$, we denote by $\sigma(T)$ the spectrum in $C^*_{\text{red}}(X)$, while $\sigma_{\max}(T)$ denotes the, possibly larger, spectrum in $C^*_{\max}(X)$.
Recall that $X$ is an expander sequence if and only if there is $h>0$ such that $\sigma(\Delta_X)\subseteq \{0\}\cup [h,\infty)$.
\begin{definition}[\cite{WilYu12}]
Let $X=(X_n)$ be a sequence of finite connected graphs with uniformly bounded degree and number of vertices going to infinity.
Then $X$ has \emph{geometric property (T)} if there is $\gamma > 0$ such that $\sigma_{\max}(\Delta_X) \subseteq \{0\}\cup [\gamma,\infty)$.
\end{definition}

A sequence of finite graphs $(X_n)$ has \emph{large girth} if for every $R>0$ there is a positive integer $N$ such that for $n\geq N$, the graphs $X_n$ do not have any cycles of length shorter than $R$.
It was shown in \cite{WilYu12} that a sequence of graphs with property (T) can never have large girth, using properties of K-theory.
We give a quantitative version of this result: there is some $R$ such that all of the $X_n$ have ``many'' $R$-cycles.

Let us make this more precise.
As in \cite{Elek06}, we define the cycle spaces of a graph.

\begin{definition}\label{def-cycle-spaces}
Let $X$ be a graph with edge set $E$.
Let $\mathbb C[E]$ be the free $\mathbb R$-vector space generated by $E$.
We use the convention that $(x,y)=-(y,x)$.
Let $Z(X)$ be the subspace of $\mathbb C[E]$ generated by the cycles, where a cycle $(x_1,x_2,\ldots,x_n)$ corresponds to the element $(x_1,x_2)+(x_2,x_3)+\ldots+(x_{n-1},x_n)+(x_n,x_1)$.
For any integer $R$, let $Z_R(X)$ be the subspace generated by the cycles of length at most $R$.
\end{definition}

We can now state the first theorem of this paper.

\begin{restatable}{theorem}{manyRcycles}\label{thm-many-R-cycles}
Suppose $X=(X_n)$ is a sequence of finite connected graphs with uniformly bounded degree and with number of vertices going to infinity.
Suppose that $X$ has geometric property (T).
Then there are constants $R>0$ and $\epsilon>0$ such that for all large enough $n$, we have $\dim Z_R(X_n)\geq\epsilon|V(X_n)|$.
\end{restatable}

The \emph{cost} of generating an equivalence relation is introduced by Levitt in \cite{Lev95}.
Taking the supremum of all equivalence relations generated by probability measure preserving actions of a group gives the cost of this group, as introduced by Gaboriau in \cite{Gab00}.
\emph{Combinatorial cost} is a variant on this concept for sequences of graphs which was defined by Elek in \cite{Elek06}.
It measures the number of edges necessary to induce the coarse structure of a sequence of graphs.
\begin{definition}
Let $(X_n)$ and $(Y_n)$ be sequences of finite connected graphs, such that $X_n$ and $Y_n$ have the same vertex set.
For $x,y\in V(X_n)$ denote by $d_X(x,y)$ the path distance in $X_n$, and by $d_Y(x,y)$ the path distance in $Y_n$.
We say that $(X_n)$ and $(Y_n)$ \emph{induce the same coarse structure} if there is a constant $L$ such that $d_X(x,y)\leq Ld_Y(x,y)$ and $d_Y(x,y)\leq Ld_X(x,y)$ for all $n$ and $x,y\in V(X_n)$.
\end{definition}

\begin{definition}[\cite{Elek06}]\label{def-combinatorial-cost}
The \emph{combinatorial cost} of a sequence $X=(X_n)$ of graphs is defined as
\[c(X) = \inf_Y \liminf_n\frac{|E(Y_n)|}{|V(X_n)|}\]
where the infimum is taken over all sequences of graphs $(Y_n)$ on the same vertex sets as $X_n$, which induce the same coarse structure as $(X_n)$.
\end{definition}

In \cite{HutchPete20}, it was shown that any group with property (T) has cost 1.
This raises the question if a similar theorem can be proved in the combinatorial context.
From Theorem \ref{thm-many-R-cycles} it already follows that the infimum in Definition \ref{def-combinatorial-cost} is not attained for sequences of graphs with geometric property (T) (see Proposition \ref{prop-cost-not-attained}).
The question, whether sequences of graphs with geometric property (T) necessarily have cost 1, remains open.

Theorem \ref{thm-many-R-cycles} also raises the following question: is it true that for every sequence $(X_n)$ with geometric property (T), there is an $R$ such that all vertices have a cycle in their $R$-neighbourhood?
We show that if two graph sequences $(X_n)$ and $(Y_n)$ ``mostly'' agree, and one of them has geometric property (T) and the other one is an expander, then the other one has geometric property (T) too.
This can be used to answer the above question negatively (see Corollary \ref{cor-geometric-T-with-locally-large-girth}).

\begin{definition}\label{def-approximating-graphs}
Let $X=(X_n)$ and $Y=(Y_n)$ be sequences of finite connected graphs with increasing number of vertices and uniformly bounded degree.
We say that $X$ and $Y$ are \emph{approximately isomorphic} if there are subgraphs $X'_n\subseteq X_n$ and $Y'_n\subseteq Y_n$ such that $X'_n$ and $Y'_n$ are isomorphic for all $n$, and
\[\lim_{n\to\infty}\frac{|V(X'_n)|}{|V(X_n)|} = \lim_{n\to\infty}\frac{|E(X'_n)|}{|E(X_n)|} = \lim_{n\to\infty}\frac{|V(Y'_n)|}{|V(Y_n)|} = \lim_{n\to\infty}\frac{|E(Y'_n)|}{|E(Y_n)|}=1.\]
\end{definition}
It is easy to see that we may take $X'_n$ and $Y'_n$ to be induced subgraphs (that is, subgraphs that contain all edges between their vertices that are present in the ambient graph).

If $X$ and $Y$ are approximately isomorphic and $X$ has geometric property (T), then $Y$ need not have geometric property (T); indeed, it does not even need to be an expander (it is only an asymptotic expander, as defined in \cite{KhuLiVigZha21}).
However, this is the only thing that can go wrong.

\begin{restatable}{theorem}{approximatingpropT}\label{thm-approximating-prop-T}
Let $X=(X_n)$ and $Y=(Y_n)$ be approximately isomorphic sequences of finite connected graphs with uniformly bounded degree and number of vertices going to infinity.
Suppose $X$ has geometric property (T) and $Y$ is an expander.
Then $Y$ has geometric property (T).
\end{restatable}

\section{Short-cycle spaces in graphs with geometric property (T)}

Let $(X_n)$ be a sequence of graphs.
The algebraic Roe algebra $\mathbb C_{\cs}[X]$, as recalled in Section 1, has a standard representation on $L^2X$.
Other representations are harder to describe explicitly.
We will see how to construct a representation of $\mathbb C_{\cs}[X]$, starting from a representation of only part of the matrix algebra $B(L^2X_n)$, for each $n$.

\begin{definition}
Let $X$ be a graph and let $B_R(L^2X)$ be the set of bounded operators on $L^2X$ of propagation at most $R$.
An $R$-representation is a Hilbert space $\H$ together with a linear map $\pi\colon B_R(L^2X)\to B(\H)$ satisfying $\pi(T^*)=\pi(T)^*$ and $\pi(TS)=\pi(T)\pi(S)$ if $T$, $S$ and $TS$ have propagation at most $R$.
\end{definition}

Suppose $X=(X_n)$ is a sequence of graphs and suppose we have a sequence of $r_n$-representations $\pi_n$ of $X_n$ on $\H_n$.
If $r_n\to\infty$, we can make a representation of $\mathbb C_{\cs}[X]$ as follows: let $\mathcal U$ be a free ultrafilter on $\mathbb N$ and define the ultraproduct $\H = \prod_n\H_n/\mathcal U$.
Define $\pi(T)=\lim_\mathcal U \pi_n(T|_{L^2X_n}) \in B(\H)$.
Then $\pi\colon\mathbb C_{\cs}[X]\to B(\H)$ is a representation.

We can use these representations to prove that a sequence of graphs with property (T) does not have large girth (which was also already proved in \cite{WilYu12} using K-theory).
Below we first give the proof below for $d$-regular graphs with an even integer $d$, to avoid technicalities.
It is true in greater generality; this will follow later from Theorem \ref{thm-many-R-cycles}.
\begin{proposition}
Let $X=(X_n)$ be a sequence of $d$-regular finite connected graphs and an increasing number of vertices.
Suppose that $X$ has geometric property (T) and that $d$ is even.
Then $X$ does not have large girth.
\end{proposition}
\begin{proof}
Suppose that $X$ has large girth and let $t\in\mathbb R$.
Let $r_n$ denote the girth of $X_n$.
There is an Eulerian cycle on $X_n$.
This defines a direction on each edge of $X_n$.
So we can define a function $\rho\colon E_n\to \{-1,1\}$, where $E_n$ denotes the edge set of $X_n$, such that for each $x\in X_n$ we have $\sum_{y\sim x}\rho(x,y)=0$.
For all pairs $(x,y)\in X_n^2$ choose a shortest path $\gamma_{(x,y)}=(x=x_0,x_1,\ldots,x_{d(x,y)}=y)$ from $x$ to $y$.
Choose it in such a way that $\gamma_{(y,x)}$ is the inverse of $\gamma_{(x,y)}$.
Define $\rho\colon X_n^2\to\mathbb Z$ by $\rho(x,y)=\rho(x_0,x_1)+\ldots+\rho(x_{d(x,y)-1},x_{d(x,y)})$.
Then $\rho(x,y)=-\rho(y,x)$, and if $d(x,y)+d(y,z)+d(x,z)<r_n$, we have $\rho(x,y)+\rho(y,z)=\rho(x,z)$.
Now define $\pi_n\colon B(L^2X_n)\to B(L^2X_n)$ by 
\[\pi_n(T)\xi(x)=\sum_yT_{xy}\exp(it\rho(x,y))\xi(y).\]
Then $\pi_n(T^*)=\pi_n(T)^*$, and if the propagations of $S$ and $T$ and $TS$ are at most $\frac13r_n$, then $\pi_n(TS)=\pi_n(T)\pi_n(S)$.
So $\pi_n$ defines a $\frac13r_n$-representation.

Consider the constant unit vector $\xi_n\in L^2X_n$ with $\xi_n(x)=\frac1{\sqrt{|V(X_n)|}}$.
This satisfies $\pi_n(\Delta_X)\xi_n=d(1-\cos t)\xi_n$.
Let $\pi:\mathbb C_{\cs}[X]\to B(\H)$ denote the limit of the $\pi_n$.
This is a representation.
Let $\xi=(\xi_n)\in \H$.
Then $\pi(\Delta_X)\xi=d(1-\cos t)\xi$, so $d(1-\cos t)\in \sigma_{\max}(\Delta_X)$.
Since $t$ can be any real number, we conclude that $[0,2d]\subseteq \sigma_{\max}(\Delta_X)$.
Then $X$ does not have geometric property (T), and we have arrived at a contradiction.
\end{proof}

If the graph is not $d$-regular for an even $d$, we do not necessarily have a Eulerian cycle.
In this case, we can still prove the result.
We need Lemmas \ref{lem-finding-a-nice-cycle} and \ref{lem-not-too-many-bridges} below to construct a cycle $v\in Z(X)$ that will play a similar role as the Eulerian cycle above.
We will use it to prove something even stronger, namely that there are ``many'' small cycles in $X$ (Theorem \ref{thm-many-R-cycles}).

Note that, if $X$ is a finite connected, then $\dim(\mathbb C[E])=|E|$ and $\dim(Z(X))=|E|-|X|+1$ (this is clear for a tree, and any time an edge is added to a graph, both sides increase by 1).
We will also use the standard inner product on $\mathbb C[E]$, such that edges have norm 1 and are perpendicular to each other.
An edge in a connected graph is called a \emph{bridge} if the removal of this edge would render the graph disconnected.

\begin{lemma}\label{lem-finding-a-nice-cycle}
Let $X$ be a finite connected graph with edge set $E$.
There is $v\in Z(X)$ satisfying the following conditions:
\begin{enumerate}[(i)]
    \item For all $e\in E$ we have $v(e)\in\{-1,0,1\}$.
    \item Of all edges $e\in E$ that are not bridges, at most half satisfy $v(e)=0$.
\end{enumerate}
\end{lemma}
\begin{proof}
Let $Z_1,\ldots,Z_K$ be cycles in $X$ such that each edge that is not a bridge is contained in one of the cycles $Z_i$.
Choose $\epsilon_1,\ldots,\epsilon_K\in\{0,1\}$ uniformly and independently at random.
Let $w=\sum_{i=1}^K\epsilon_iZ_i$.
For each $e\in E$ that is not a bridge, we have $\mathbb P[w(e)\text{ is odd}]=\frac12$.
So $\mathbb E[\#\{e\in E\mid w(e)\text{ is odd}\}]$ is equal to half the number of edges that are not bridges.
So there is some $w\in Z(X)$ such that $w(e)$ is odd for at least half the number of edges that are not bridges.
Now let $E'\subseteq E$ be the subset of $E$ consisting of the edges for which $w(e)$ is odd.
Then each vertex of $X$ has an even number of adjacent edges in $E'$.
So each component of $(V(X),E')$ has a Eulerian cycle.
Let $v$ be the sum of these Eulerian cycles.
Then $v$ satisfies the conditions.
\end{proof}

\begin{lemma}\label{lem-not-too-many-bridges}
Let $X$ be a finite connected graph with edge set $E$ and maximal degree $d$.
Suppose there is a constant $h>0$ such that for all subsets $A$ with $\frac14|V(X)|\leq |A|\leq\frac12|V(X)|$, we have $|\delta A|\geq h|A|$.
Then there is a constant $c>0$, only depending on $d$ and $h$, such that if $|V(X)|$ is large enough, the number of edges that are not bridges in $X$ is at least $c|E|$.
\end{lemma}
\begin{proof}
For $b$ a bridge, define $K_b$ and $G_b$ to be the two components of the graph $(V(X),E\setminus\{b\})$, with $|K_b|\leq|G_b|$.
Since $\delta K_b=\{b\}$, and $h\cdot \frac14|V(X)|>1$ for $|V(X)|$ large enough, we must have $|K_b|<\frac14|V(X)|$.
Now let $G=\cap_bG_b$ and $K=V(X)\setminus G$, where the intersection ranges over all bridges $b$.

We show that $G$ is non-empty.
Let $b$ be a bridge such that $K_b$ is maximal.
Let $x$ be the endpoint of $b$ in $G_b$.
Then $x\in G$: for if there is another bridge $b'$ with $x\in K_{b'}$, we have either $K_b\cup\{x\}\subseteq K_{b'}$, or $G_b\subseteq K_{b'}$, and then $K_{b'}$ is too large.

Note that all bridges of $X$ must have at least one vertex in $K$.
Therefore the number of edges that are not bridges is at least $|G|-1$.
So if we find a constant $c'>0$ such that $|G|\geq c'|V(X)|$, we are done.

Let $\delta G=\{b_1,\ldots,b_N\}$.
These are all bridges.
For each $1\leq i\leq N$ we have $|K_{b_i}|\leq\frac14|V(X)|$.
Then we can choose some $1\leq M\leq N$ such that $\min(|K|,\frac14|V(X)|)\leq\sum_{j=1}^M|K_{b_j}|\leq\frac12|V(X)|$.
If $K<\frac14|V(X)|$, we have $|G|\geq\frac34|V(X)|$, and we are done.
Suppose that $K\geq\frac14|V(X)|$.
Then we can apply the assumption of the lemma to $A=\cup_{j=1}^MK_{b_j}$.
We have $|\delta A|=M\leq N=|\delta G|\leq d|G|$, therefore $|A|\leq \frac dh|G|$.
So $|G|\geq \frac h{4d}|V(X)|$.
This finishes the proof.
\end{proof}

\begin{corollary}\label{cor-finding-a-nice-cycle}
Let $X$ be a finite connected graph with edge set $E$ and maximal degree $d$.
Suppose there is a constant $h>0$ such that for all subsets $A\subset V(X)$ with $\frac14|V(X)|\leq |A|\leq\frac12|V(X)|$, we have $|\delta A|\geq h|A|$.
Then, provided that $|V(X)|$ is large enough, there are a constant $c>0$, only depending on $d$ and $h$, and $v\in Z(X)$ satisfying the following conditions:
\begin{enumerate}[(i)]
    \item For all $e\in E$ we have $v(e)\in\{-1,0,1\}$.
    \item We have $\#\{e\in E\mid v(e)\neq 0\}\geq c|E|$.
\end{enumerate}
\end{corollary}
\begin{proof}
This follows directly from Lemmas \ref{lem-finding-a-nice-cycle} and \ref{lem-not-too-many-bridges}.
\end{proof}

\manyRcycles*
\begin{proof}
Since $X$ is an expander sequence, it can have only finitely many trees.
We can then assume without loss of generality that none of the $X_n$ is a tree.

Let $\gamma > 0$ be such that $\sigma_{\max}(\Delta_X) \subseteq \{0\}\cup[\gamma,\infty)$.
Let $d$ be the maximal degree of all vertices of $X$.
Since $X$ has geometric property (T), it is in particular an expander sequence: for each subset $A\subseteq X_n$ with $|A|\leq\frac12|V(X_n)|$ we have $|\delta A|\geq\frac{\gamma^2}{4d}|A|$.

Define $h=\frac{\gamma^2}{8d}$.
Let $c_1$ be the constant from Corollary \ref{cor-finding-a-nice-cycle}, using the constants $d$ and $h$.
Let $c_2>0$ be such that
\begin{equation}\label{ineq-with-c_2-and-c_3}c_3=\frac14d^2-\left(1-\frac{c_1}{2d}\right)\left(\frac12d+c_2d\right)^2-\frac{c_1}{2d}\left(\frac12d+c_2d-\frac12\right)^2>0.\end{equation}
Let $t>0$ be small enough such that
\begin{equation}\label{ineq-important-with-t}dt^2<\gamma\end{equation}
and also
\begin{equation}\label{ineq-exponential-with-t}\left|\exp(it)-1-it+\frac12t^2\right|\leq c_2t^2.\end{equation}
Let $\epsilon>0$ be such that the following conditions are satisfied:
\begin{align}
8\epsilon d^2&\leq\frac12c_3t^4,\label{ineq-with-epsilon}\\
2\epsilon&\leq\frac{c_1}{2d},\label{ineq-also-with-epsilon}\\
4\epsilon&\leq h\label{ineq-another-one-with-epsilon}.
\end{align}

Suppose for a contradiction that for each $R$, there is an $n$ such that $\dim Z_R(X_n)<\epsilon|V(X_n)|$.
Let $E_n$ be the set of edges of $X_n$.
Note that $Z(X_n)\subseteq \mathbb C[E_n]$ consists of the functions $\rho:E_n\to \mathbb C$ satisfying $\sum_{y\sim x}\rho(x,y)=0$ for all $x\in V(X_n)$.
The subset $Z_R(X_n)^\perp\subseteq \mathbb C[E_n]$ consists of all function $\rho\colon E_n\to \mathbb C$ satisfying $\rho(x_1,x_2)+\rho(x_2,x_3)+\ldots+\rho(x_q,x_1)=0$ for all $q$-cycles $(x_1,x_2,\ldots,x_q)$ with $q\leq R$.

Let $\rho\in Z_R(X_n)^\perp$.
We will construct a $\frac13R$-representation $\pi_\rho$ of $X_n$.
First, we extend $\rho$ to a function on $V(X_n)^2$ as follows: for each pair $(x,y)\in V(X_n)^2$, choose a shortest path $(x=x_0,x_1,\ldots,x_{d(x,y)}=y)$ from $x$ to $y$ and define $\rho(x,y)=\rho(x_0,x_1)+\rho(x_1,x_2)+\ldots+\rho(x_{d(x,y)-1},x_{d(x,y)})$.
Since $\rho\in Z_R(X_n)^\perp$, this satisfies the following: if $x,y\in V(X_n)$ satisfy $d(x,y)\leq R$, then $\rho(x,y)=-\rho(y,x)$, and if $x,y,z\in V(X_n)$ satisfy $d(x,y)+d(y,z)+d(z,x)\leq R$, then $\rho(x,y)+\rho(y,z)+\rho(z,x)=0$.

Now we define $\pi_\rho:B(L^2X_n)\to B(L^2X_n)$ by
\[\pi_\rho(T)\xi(x)=\sum_yT_{xy}\exp(i\rho(x,y))\xi(y).\]
This is a $\frac13R$-representation.

For a subset $B\subseteq E_n$, denote by $X_n\setminus B$ the graph $(X_n,E_n\setminus B)$.
Choose edges $b_1,\ldots,b_{\dim(Z_R(X_n))}$ recursively in such a way that $b_j$ is an edge in a cycle in $Z_R(X_n)\cap Z(X_n\setminus\{b_1,\ldots,b_{j-1}\})$.
Then $Z_R(X_n)\cap Z(X_n\setminus\{b_1,\ldots,b_{j-1},b_j\})$ has one dimension fewer than $Z_R(X_n)\cap Z(X_n\setminus\{b_1,\ldots,b_{j-1}\})$.
Thus, with $B=\{b_1,\ldots,b_{\dim(Z_R(X_n))}\}$ we have $Z_R(X_n)\cap Z(X_n\setminus B)=0$.
Then we also have $Z_R(X_n)\cap \mathbb C[E_n\setminus B]=0$, and counting dimensions, $Z_R(X_n)\oplus \mathbb C[E_n\setminus B]=\mathbb C[E_n]$.
Since $\mathbb C[B]\perp\mathbb C(E_n\setminus B)$ and $Z_R(X_n)^\perp\perp Z_R(X_n)$, it follows that $\mathbb C[B]\cap Z_R(X_n)^\perp=0$.
Counting dimensions again, we get $Z_R(X_n)^\perp\oplus\mathbb C[B]=\mathbb C[E_n]$.

Since each $b_j$ is an edge in a cycle in $Z(X_n\setminus\{b_1,\ldots,b_{j-1}\})$, it follows by induction that $X_n\setminus B$ is still a connected graph.
Let $A\subset V(X_n)$, with $\frac14|V(X_n)|\leq |A|\leq\frac12|V(X_n)|$.
Denote by $\delta_{X_n\setminus B}A$ the set of edges in $X_n\setminus B$ with exactly one vertex in $A$.
Then we have $|\delta_{X_n\setminus B}|\geq|\delta_{X_n}A|-|B|\geq\frac{\gamma^2}{4d}|A|-\epsilon|V(X_n)|\geq(\frac{\gamma^2}{4d}-4\epsilon)|A|\geq h|A|$ by inequality (\ref{ineq-another-one-with-epsilon}).
So $X_n\setminus B$ satisfies the conditions of Corollary \ref{cor-finding-a-nice-cycle}.
We can also assume $V(X_n)$ is large enough to apply the corollary, by taking $R$ large enough.

Now let $v\in C(X_n\setminus B)$ be as in Corollary \ref{cor-finding-a-nice-cycle}.
Recall that $Z_R(X_n)^\perp\oplus\mathbb C[B]=\mathbb C[E_n]$.
Hence there is a unique function $\rho\in Z_R(X_n)^\perp$ such that $\rho(e)=v(e)$ for all $e\not\in B$.

Now we consider the $\frac13R$-representation $\pi_{t\rho}$.
Let $\Delta_t=\pi_{t\rho}(\Delta_X)\in B(L^2X_n)$.
Let $\xi_n\in L^2X_n$ be the constant function with $\xi_n(x)=1$ for all $x\in X_n$.
Then we have
\[\Delta_t\xi_n(x)=\sum_{y\sim x}(1-\exp(it\rho(x,y)))\]
for all $x\in X_n$.

We define the following subsets of $V(X_n)$:
\begin{align*}
    A_3&=\{x\in X_n\mid (x,y)\in B\text{ for some }y\sim x\},\\
    A_1&=\{x\in X_n\setminus A_3\mid v(x,y)\neq 0\text{ for some }y\sim x\},\\
    A_2&=X_n\setminus(A_1\cup A_3).
\end{align*}
Since $|B|\leq\epsilon|V(X_n)|$, we have $|A_3|\leq2\epsilon|V(X_n)|$, and since at least $c_1|E_n|$ edges $e$ satisfy $v(e)\neq0$, we have $|A_1|\geq\frac{c_1}d|E_n|-|A_3|\geq(\frac{c_1}d-2\epsilon)|V(X_n)|\geq\frac{c_1}{2d}|V(X_n)|$ by inequality (\ref{ineq-also-with-epsilon}).

We will estimate $|\Delta_t\xi_n(x)-\frac12dt^2|$, finding increasingly good estimates in the cases $x\in A_3,A_2,A_1$ respectively.
For all $x\in X_n$, we have $|\Delta_t\xi_n(x)-\frac12dt^2|\leq 2d$.
For $x\not\in A_3$ we have $\rho(x,y)=v(x,y)$ for all $y\sim x$.
Since $\sum_{y\sim x}v(x,y)=0$, we have
\[\Delta_t\xi_n(x)=\sum_{y\sim x}(1+itv(x,y)-\exp(itv(x,y))).\]
Then, using inequality (\ref{ineq-exponential-with-t}) we have:
\begin{align*}
    \left|\Delta_t\xi_n(x)-\frac12dt^2\right|&\leq\sum_{y\sim x}\left|1+itv(x,y)-\frac d{2\deg(x)}t^2-\exp(itv(x,y))\right|\\
    &\leq\left(\sum_{y\sim x}\left|1+itv(x,y)-\frac d{2\deg(x)}t^2-(1+itv(x,y)-\frac12t^2v(x,y)^2)\right|\right)+c_2dt^2\\
    &\leq \frac12dt^2+c_2dt^2-\frac12\sum_{y\sim x}v(x,y)^2t^2.
\end{align*}
For $x\in A_2$, this is at most $(\frac12d+c_2d)t^2$.
For $x\in A_1$, there is $y\sim x$ with $v(x,y)\neq 0$, and we find that $|\Delta_t\xi_n(x)-\frac12dt^2|\leq(\frac12d+c_2d-\frac12)t^2$.

Using inequalities (\ref{ineq-with-c_2-and-c_3}) and (\ref{ineq-with-epsilon}), we now have:
\begin{align*}
    \norm{\Delta_t\xi_n-\frac12dt^2\xi_n}^2&=\sum_{x\in X_n}|\Delta_t\xi_n(x)-\frac12dt^2|^2\\
    &\leq\sum_{x\in A_1}\left(\frac12d+c_2d-\frac12\right)^2t^4+\sum_{x\in A_2}\left(\frac12d+c_2d\right)^2t^4+\sum_{x\in A_3}4d^2\\
    &\leq\left(\left(1-\frac{c_1}{2d}\right)\left(\frac12d+c_2d\right)^2+\frac{c_1}{2d}\left(\frac12d+c_2d-\frac12\right)^2\right)|V(X_n)|t^4+8\epsilon d^2|V(X_n)|\\
    &=\left(\frac14d^2-c_3\right)|V(X_n)|t^4+8\epsilon d^2|V(X_n)|\\
    &\leq \left(\frac14d^2-\frac12c_3\right)|V(X_n)|t^4\\
    &=\left(1-\frac{2c_3}{d^2}\right)\norm{\frac12dt^2\xi_n}^2.
\end{align*}

For each $R$ we have found an $n_R$, a $\frac13R$-representation $\pi_{n_R}:B(L^2X_{n_R})\to B(L^2X_{n_R})$, and a vector $\xi_{n_R}\in B(L^2X_{n_R})$ satisfying $\norm{\pi_{n_R}(\Delta_X)\xi_{n_R}-\frac12dt^2\xi_{n_R}}\leq(1-\frac{2c_3}{d^2})^\frac12\norm{\frac12dt^2\xi_{n_R}}$.
Let $\mathcal U$ be a free ultrafilter on $\mathbb N$ and define the ultraproduct $\H=\prod_nB(L^2X_{n_R})/\mathcal U$.
Let $\pi=\lim_\mathcal U\pi_{n_R}$ and $\xi=\lim_\mathcal U\frac{\xi_{n_R}}{\norm{\xi_{n_R}}}$.
Then $\pi$ is a representation of $\mathbb C_{\cs}[X]$, and we have $\norm{\pi(\Delta_X)\xi-\frac12dt^2\xi}\leq(1-\frac{2c_3}{d^2})^\frac12\norm{\frac12dt^2\xi}$.
It follows that $\sigma(\pi(\Delta_X))\cap[\frac12dt^2(1-(1-\frac{2c_3}{d^2})^\frac12),\frac12dt^2(1+(1-\frac{2c_3}{d^2})^\frac12)]\neq\emptyset$.
So $\sigma_{\max}(\Delta_X)$ contains a positive element that is at most $dt^2$.
This is a contradiction with inequality (\ref{ineq-important-with-t}).
\end{proof}

\begin{remark}

Note that in the above proof, $\epsilon$ only depends on $d$ and $\gamma$.
\end{remark}

We have shown that each sequence of graphs with geometric property (T) has many small cycles.
It follows that we can remove a large number of cycles of the graph, while still keeping the same coarse structure.
In particular, the infimum in the definition of cost is not attained (see Definition \ref{def-combinatorial-cost}).

\begin{proposition}\label{prop-cost-not-attained}
Let $X=(X_n)$ be a sequence of graphs with degree at most $d$.
Suppose there are $R,\epsilon>0$ such that $\dim Z_R(X_n)\geq \epsilon|V(X_n)|$ for all $n$.
Then $c(X)\leq \liminf_n\frac{|E(X_n)|}{|V(X_n)|}-\frac{\epsilon}{d^{R-1}}$.
In particular, if $X$ has geometric property (T), the infimum $c(X)=\inf_{Y}\liminf_n\frac{|E(Y_n)|}{|V(X_n)|}$, over all sequences $(Y_n)$ inducing the same coarse structure as $X$, is not attained.
\end{proposition}
\begin{proof}
Consider the set of all subgraphs of $X_n$, on the same vertex set, without $R$-cycles.
Let $Y_n$ be a maximal element of this set.
If $x,y\in X_n$ are connected, then either they are connected in $Y_n$, or adding the edge $(x,y)$ to $Y_n$ would create an $R$-cycle in $Y_n$.
So $d_Y(x,y)\leq R-1$.
This shows that $X$ and $Y$ induce the same coarse structure.

Since $\dim Z_R(X_n)\geq \epsilon|V(X_n)|$, there are in particular at least $\epsilon |V(X_n)|$ cycles of length at most $R$ in $X_n$.
Of all these cycles, at least one of its edges is not in $Y_n$.
One edge can be contained in at most $d^{R-1}$ cycles of length at most $R$.
So $|E(Y_n)|\leq |E(X_n)|-\frac{\epsilon}{d^{R-1}}|V(X_n)|$.
This gives $\liminf_n\frac{|E(Y_n)|}{|V(X_n)|}\leq \liminf_n\frac{|E(X_n)|}{|V(X_n)|} - \frac{\epsilon}{d^{R-1}}$.
So $c(X)\leq \liminf_n\frac{|E(X_n)|}{|V(X_n)|}-\frac{\epsilon}{d^{R-1}}$.

The last statement of the lemma follows from Theorem \ref{thm-many-R-cycles} and the fact that if $X$ and $Y$ induce the same coarse structure, also $Y$ has geometric property (T) (see \cite[Theorem 4.3]{WilYu14}).
\end{proof}

\section{Behaviour of geometric property (T) under small changes}
A natural question in light of Theorem \ref{thm-many-R-cycles} is the following: if a sequence of graphs $(X_n)$ has geometric property (T), is it necessary that there is an $R$ such that every vertex has a cycle in its $R$-neighbourhood?
In this section we will see that the answer is no.
We will do this by proving that if we change a sequence of graphs a small amount (see Definition \ref{def-approximating-graphs}) while keeping expansion, we also keep geometric property (T).

\approximatingpropT*
For the proof we need the following proposition.
\begin{proposition}\label{prop-almost-constant-vector-in-expander}
Let $X=(X_n)$ be an expander sequence with maximum degree $d$, and let $h>0$ such that $\sigma(\Delta_X)\subseteq \{0\}\cup[h,\infty)$.
Let $\pi\colon \mathbb C_{\cs}[X]\to B(\H)$ be a representation and suppose that $v\in\H$ is a unit vector with $\Delta_Xv=\eta v$ for some $\eta>0$.
\begin{enumerate}[(i)]
    \item Let $F_n\subseteq V(X_n)$ be such that $\lim_{n\to\infty}\frac{|F_n|}{|V(X_n)|} = 0$, and let $P_F\in \mathbb C_{\cs}[X]$ denote the projection on the vertices of the $F_n$.
Then $\norm{P_Fv}\leq 2^\frac34d^\frac12h^{-\frac12}\eta^\frac14$.
\item Let $G=(G_n)$ be a sequence of subgraphs of $X_n$ such that $\lim_{n\to\infty}\frac{|V(G_n)|}{|V(X_n)|}=0$.
Then $\langle\Delta_Gv,v\rangle\leq 2^\frac32dh^{-1}\eta^\frac32$ and $\norm{\Delta_Gv}\leq 2^\frac34d^\frac12h^{-\frac12}\eta^\frac54$.
\end{enumerate}
\end{proposition}
\begin{proof}
\begin{enumerate}[(i)]
\item We can assume that $\eta<h$.
Let $\delta>0$ and let $N$ be an integer such that $|F_n|/|V(X_n)|<\delta $ for $n>N$.
Let $P_{X_{\leq N}}\in \mathbb C_{\cs}[X]$ be the projection on the vertices of $X_1$ up to $X_N$.
Then we have $\Delta_XP_{X_{\leq N}}v = \eta P_{X_{\leq N}}v$.
Since $\eta$ is not an eigenvalue of $\Delta_XP_{X_{\leq N}}$, we have $P_{X_{\leq N}}v=0$.
So $P_Fv = P_{F_{>N}}v$, where $P_{F_{>N}}$ is the projection on the union $F_{>N}=\bigcup_{n>N}F_n$.

We can colour the edges of $X$ in $2d$ colours such that no two intersecting edges have the same colour.
For each colour $i$ define the involution $\tau_i$ of $X$ that sends each vertex in an edge with colour $i$ to the other vertex of this edge, while fixing the other vertices.
This defines an element in $\mathbb C_{\cs}[X]$ that we also denote by $\tau_i$.
Since $\tau_i^2=1$, we have $0\leq 1-\tau_i\leq 2$.
This inequality holds in $C^*_{\max}(X)$.
We have $\Delta_X = \sum_{i=1}^{2d}(1-\tau_i)$.

For $f\in L^\infty X$, denote the corresponding multiplication operator in $B(L^2X)$ by $M_f$.
Define the positive unital linear map $\phi\colon L^\infty X\to \mathbb C$ by $\phi(f) = \langle M_fv,v\rangle$.
For each $i$ we have $\langle(1-\tau_i)v,v\rangle \leq \eta$, so $\norm{(1-\tau_i)v}^2=\langle(2-2\tau_i)v,v\rangle\leq 2\eta$.
Then for $f\in L^\infty X$ we have
\begin{align*}
    (1-\tau_i)\phi(f)& = \phi(f-f\circ\tau_i)\\
    &=\langle M_fv,v\rangle - \langle M_{f\circ\tau_i}v,v\rangle\\
    &=\langle M_fv,v\rangle - \langle M_f\tau_iv,\tau_iv\rangle\\
    &=\langle M_fv,v-\tau_iv\rangle + \langle M_f(v-\tau_iv),\tau_iv\rangle.
\end{align*}
We get $|(1-\tau_i)\phi(f)|\leq 2\norm{f}_\infty\norm{v-\tau_iv}\leq 2\sqrt{2\eta}\norm{f}_\infty$, so $\norm{(1-\tau_i)\phi}\leq 2\sqrt{2\eta}$.

Let $\mathcal P(X)\subseteq L^1X$ denote the set of probability measures on $X$, and let $\mathcal B$ be the real Banach space $\bigoplus_{i=1}^{2d}L^1(X,\mathbb R)$ with norm $\norm{(\psi_i)} = \max_i\norm{\psi_i}_1$.
Define the convex set $B\subseteq \mathcal B$ by
\[B = \{((1-\tau_i)\psi) \mid \psi\in\mathcal P(X), |\psi(\mathbb 1_{F_{>N}})-\phi(\mathbb 1_{F_{>N}})| < \delta\}.\]
Let $A$ be the open ball $\{ a \in \mathcal B \mid \norm{a}<2\sqrt{2\eta}+\delta\}$.
Suppose $A$ and $B$ are disjoint.
By the Hahn-Banach separation theorem, there is $f\in B^*$ and a positive real number $s$ with $f(a)<s\leq f(b)$ for $a\in A$ and $b\in B$.
By the Goldstine theorem, $\phi$ is in the weak closure of $\mathcal P(X)$.
Then $((1-\tau_i)\phi)\in B^{**}$ is in the weak closure of $B$, so $f(((1-\tau_i)\phi))\geq s$.
On the other hand, $\norm{((1-\tau_i)\phi)}\leq 2\sqrt{2\eta}$, so $((1-\tau_i)\phi)$ is in the weak closure of $\frac{2\sqrt{2\eta}}{2\sqrt{2\eta}+\delta}A$, showing that $f(((1-\tau_i)\phi))\leq \frac{2\sqrt{2\eta}}{2\sqrt{2\eta}+\delta}s$, giving a contradiction.

Therefore, there is $\psi\in A\cap B$.
Let $\xi = \psi^\frac12 \in L^2X$.
Then for all $i$, we have
\begin{align*}\langle (1-\tau_i)\xi,\xi\rangle &= \sum_{x\in X}(\xi(x)-\xi(\tau_i x))\xi(x)\\
&=\frac12\sum_{x\in X}\left(\xi(x)-\xi(\tau_i x)\right)^2\\
&\leq \frac12\sum_{x\in X}|\xi(x)^2-\xi(\tau_i x)^2|\\
&=\frac12\norm{(1-\tau_i)\psi}_1\\
&\leq \sqrt{2\eta}+\frac12\delta.
\end{align*}
Summing over all $i$ gives $\langle \Delta_X\xi,\xi\rangle \leq 2d\sqrt{2\eta}+d\delta$.
Let $\xi_c$ be the projection of $\xi$ on the locally constant functions in $L^2X$.
We get $\langle \Delta_X(\xi-\xi_c),\xi-\xi_c\rangle = \langle \Delta_X\xi,\xi\rangle\leq 2d\sqrt{2\eta}+d\delta$, but also $\langle\Delta_X(\xi-\xi_c),\xi-\xi_c\rangle \geq h\norm{\xi-\xi_c}_2^2$.
So $\norm{\xi-\xi_c}_2^2\leq \frac1h(2d\sqrt{2\eta}+d\delta)$.

Since $\frac{|F_n|}{V(X_n)} < \delta$ for $n>N$, we have $\norm{\xi_{c|F_{>N}}} < \sqrt\delta$.
We get
\begin{align*}
    \psi(\mathbb 1_{F_{>N}})&=\sum_{x\in F_{>N}}\xi(x)^2\\
    &=\sum_{x\in F_{>N}}(\xi-\xi_c)(x)^2 + 2\sum_{x\in F_{>N}}\xi_c(x)(\xi-\xi_c)(x) + \sum_{x\in F_{>N}}\xi_c(x)^2\\
    &\leq \norm{\xi-\xi_c}_2^2+2\norm{\xi_{c|F_{>N}}}_2\cdot\norm{\xi-\xi_c}_2 + \norm{\xi_{c|F_{>N}}}_2^2\\
    &\leq \frac1h(2d\sqrt{2\eta}+d\delta) + 2\sqrt{\delta} +\delta.
\end{align*}

Finally, we have
\[\norm{P_Fv}^2 =\norm{P_{F>N}v}^2 =\langle M_{\mathbb 1_{F_{>N}}}v,v\rangle = \phi(\mathbb 1_{F_{>N}})\leq \psi(\mathbb 1_{F_{>N}})+\delta\leq \frac1h(2d\sqrt{2\eta}+d\delta) + 2\sqrt{\delta} + 2\delta.\]
Letting $\delta\to0$ gives the desired conclusion.

\item We recursively define $F_0 = \bigcup_nV(G_n)$ and
\[F_{k+1}=\{x\in V(X)\setminus F_{\leq k}\mid x\text{ is adjacent to a vertex in }F_k\}.\]
Here we write $F_{\leq k}$ for $F_0\cup\ldots\cup F_k$.
Let $\delta>0$.
Since the $F_k$ are all disjoint, there is $k\geq 1$ with $\norm{P_{F_{k-1}\cup F_k}v}\leq\delta$.
Since the graphs have uniformly bounded degree, we have $\lim_{n\to\infty}\frac{|F_{\leq k-1}|}{|V(X_n)|} = 0$.
By part $(i)$, we have $\norm{P_{F_{\leq k-1}}v} \leq 2^\frac34d^\frac12h^{-\frac12}\eta^\frac14$.
Denote by $\Delta_{F_{\leq k}}$ the Laplacian operator of the induced subgraph with vertex set $F_{\leq k}$.
Then we have
\[P_{F_{\leq k-1}}\Delta_X = P_{F_{\leq k-1}}\Delta_{F_{\leq k}}=\Delta_{F_{\leq k}}-P_{F_k}\Delta_{F_{\leq k}}=\Delta_{F_{\leq k}}-P_{F_k}\Delta_{F_{\leq k}}P_{F_{k-1}\cup F_k}.\]
It follows that
\begin{align*}
\langle\Delta_Gv,v\rangle&\leq \langle\Delta_{F_{\leq k}}v,v\rangle \\
&= \langle P_{F_{\leq k-1}}\Delta_Xv,v\rangle + \langle P_{F_k}\Delta_{F_{\leq k}}P_{F_{k-1}\cup F_k}v,v\rangle \\
&= \eta\langle P_{F_{\leq k-1}}v,v\rangle + \norm{P_{F_k}\Delta_{F_{\leq k}}}\cdot\norm{P_{F_{k-1}\cup F_k}v}\cdot\norm v\\
&\leq 2^\frac32dh^{-1}\eta^\frac32+2d\delta.
\end{align*}
Similarly, we have
\begin{align*}
    \norm{\Delta_Gv}&\leq \norm{F_{\leq k}v}\\
    &\leq \norm{P_{F_{\leq k-1}}\Delta_Xv} + \norm{P_{F_k}\Delta_{F_{\leq k}}P_{F_{k-1}\cup F_k}v}\\
    &\leq 2^\frac34d^\frac12h^{-\frac12}\eta^\frac54 + 2d\delta.
\end{align*}
Since these inequalities hold for all $\delta>0$, the conclusion follows.
\end{enumerate}
\end{proof}

\begin{proof}[Proof of Theorem \ref{thm-approximating-prop-T}]
Since $X$ and $Y$ are approximately isomorphic, we can identify the isomorphic subgraphs and assume there are induced subgraphs $Z_n\subseteq X_n\cap Y_n$ with $\lim_{n\to\infty}\frac{|V(Z_n)|}{|V(X_n)|}=\lim_{n\to\infty}\frac{|V(Z_n)|}{|V(Y_n)|} = 1$.
Let $d$ be the maximum degree in $X\cup Y$, let $\gamma>0$ with $\sigma_{\max}(\Delta_X)\subseteq \{0\}\cup [\gamma,\infty)$ and let $h>0$ with $\sigma(\Delta_Y)\subseteq \{0\}\cup [h,\infty)$.
Suppose that $Y$ does not have geometric property (T).
Then the maximal spectrum of $\Delta_Y$ contains arbitrarily small positive numbers.
Let $0<\eta<h$ be in the maximal spectrum of $\Delta_Y$.
Then there is a representation $\rho\colon \mathbb C_{\cs}[Y] \to B(\H)$ and a unit vector $v \in \H$ with $\Delta_Yv=\eta v$.

We will first apply Proposition \ref{prop-almost-constant-vector-in-expander} to bound $v$ on some small subsets of $Y$.
Then we will construct a representation of $\mathbb C_{\cs}[X]$ containing a vector $1\tensor v$, that we will show is almost constant.
Since $X$ has geometric property (T), it follows that $1\tensor v$ is close to some constant vector.
This will give a contradiction because $v$ is perpendicular to all constant vectors.

We have
\[\langle\Delta_{Y\setminus Z}v,v\rangle\leq 2^\frac32dh^{-1}\eta^\frac32\]
and
\[\norm{\Delta_{Y\setminus Z}v}\leq 2^\frac34d^\frac12h^{-\frac12}\eta^\frac54\]
by Proposition \ref{prop-almost-constant-vector-in-expander}(ii).

Let $F = \{z\in Z\mid z\text{ adjacent to some }x\in X\setminus Z\}$.
Since the degree of the vertices of $Y$ is uniformly bounded, we have $\lim_{n\to\infty}\frac{|F\cap Y_n|}{|Y_n|}=0$.
By Proposition \ref{prop-almost-constant-vector-in-expander}(i), we have $\norm{P_Fv}\leq 2^\frac34d^\frac12h^{-\frac12}\eta^\frac14$.

Now we construct a representation of $\mathbb C_{\cs}[X]$.
Consider the map $E\colon \mathbb C_{\cs}[X]\to \mathbb C_{\cs}[Z]$ given by $E(T) = P_ZTP_Z$.
This is a conditional expectation, meaning that for $T\in\mathbb C_{\cs}[X]$ and $S\in\mathbb C_{\cs}[Z]$, we have $E(TS)=E(T)S$ and $E(ST)=SE(T)$.
We can construct the tensor product $\H'=\mathbb C_{\cs}[X] \tensor_{\mathbb C_{\cs}[Z]} \H$, as in \cite[Theorem 1.8]{Rie74}.
We repeat the construction here.
First we consider the algebraic tensor product $\mathbb C_{\cs}[X] \odot \H$.
We equip this with a conjugate symmetric form given on simple tensors by
\[\langle T_1\odot v_1,T_2\odot v_2\rangle = \langle E(T_2^*T_1)v_1,v_2\rangle.\]
It can be shown that this is positive semi-definite (see \cite[Lemma 1.7]{Rie74} and its proof).
Define the semi-norm $\norm w = \langle w,w\rangle ^\frac12$ for $w\in\mathbb C_{\cs}[X]\odot \H$.
Let $\H'=\mathbb C_{\cs}[X]\tensor_{\mathbb C_{\cs}[Z]}\H$ be the Hilbert space we get by taking the quotient with respect to the kernel of $\norm\cdot$ and then taking the completion.
Let $T_1\tensor v_1$ denote the image of $T_1\odot v_1$ in $\H'$.
It is easy to see that for $S\in\mathbb C_{\cs}[Z]$, we have $T_1S\tensor v_1 = T_1\tensor Sv_1$.
We now have a representation $\pi\colon \mathbb C_{\cs}[X]\to B(\H')$, given on simple tensors by $\pi(T_1)(T_2\tensor v_1) = T_1T_2\tensor v_1$.

Consider the unit vector $1\tensor v\in \H'$.
We have
\begin{align*}
    \langle \Delta_X(1\tensor v),1\tensor v\rangle &= \langle \Delta_Zv,v\rangle + \langle E(\Delta_{X\setminus Z})v,v\rangle\\
    &= \langle \Delta_Zv,v\rangle + \langle E(\Delta_{X\setminus Z})P_Fv,P_Fv\rangle\\
    &\leq \langle \Delta_Yv,v\rangle + \norm{E(\Delta_{X\setminus Z})}\cdot\norm{P_Fv}^2\\
    &\leq \eta + 2d\cdot 2^\frac32dh^{-1}\eta^\frac12\\
    &\leq 8d^2h^{-1}\eta^\frac12,
\end{align*}
provided $\eta$ is small enough.
Let $w\in \H'_c$ be the projection of $1\tensor v$ on the space of constant vectors $\H'_c = \ker(\rho(\Delta_X))$.
Then we have
\[\langle \Delta_X(1\tensor v),1\tensor v\rangle = \langle \Delta_X(1\tensor v-w),1\tensor v-w\rangle\geq \gamma\norm{1\tensor v-w}^2.\]
Combining these inequalities, we get
\[\norm{1\tensor v-w}^2 \leq 8d^2h^{-1}\gamma^{-1}\eta^\frac12.\]
Finally, we have
\begin{align*}
\eta &= \langle \Delta_Yv,v\rangle\\
&= \langle \Delta_Zv,v\rangle + \langle \Delta_{Y\setminus Z}v,v\rangle\\
&= \langle \Delta_Z(1\tensor v),w\rangle + \langle \Delta_Z(1\tensor v),1\tensor v-w\rangle +\langle \Delta_{Y\setminus Z}v,v\rangle\\
&= \langle\Delta_Z(1\tensor v),1\tensor v-w\rangle + \langle \Delta_{Y\setminus Z}v,v\rangle\\
&= \langle\eta(1\tensor v),1\tensor v-w\rangle -\langle 1\tensor\Delta_{Y\setminus Z}v,1\tensor v-w\rangle +\langle\Delta_{Y\setminus Z}v,v\rangle\\
&\leq \eta\norm{1\tensor v-w}^2 + \norm{\Delta_{Y\setminus Z}v}\cdot \norm{1\tensor v-w} + 2^\frac32dh^{-1}\eta^\frac32\\
&\leq\left( 8d^2h^{-1}\gamma^{-1} + 2^\frac94d^\frac32h^{-1}\gamma^{-\frac12} + 2^\frac32dh^{-1}\right)\eta^\frac32.
\end{align*}
This gives a contradiction if $\eta$ is small enough.
Hence $Y$ must have geometric property (T).
\end{proof}

Using the theorem, we can construct a sequence of graphs with geometric property (T) such that the graphs locally have arbitrarily large girth.

\begin{corollary}\label{cor-geometric-T-with-locally-large-girth}
There is a sequence of graphs $(X_n)$ with uniformly bounded degree and number of vertices converging to infinity, satisfying geometric property (T), with a designated vertex $p_n \in V(X_n)$, such that the $n$-ball around $p_n$ does not contain any cycles.
\end{corollary}
\begin{proof}
We start with a sequence of finite connected graphs $(Y_n)$ with geometric property (T) and maximal degree $d$.
Let $(R_n)$ be an unbounded sequence of integers with $\lim_{n\to\infty}\frac{2^{R_n}}{|V(Y_n)|} = 0$.
Let $(T_n,p_n)$ be a rooted tree of depth $R_n$ such that each vertex except for the leaves has degree 3.
Connect each leaf of the tree to a different vertex in $Y_n$.
We call the new graph $X_n$, and we show that $(X_n)$ is still an expander sequence.
Note that the maximal degree of $X_n$ equals $d+1$.
For a subset $A\subseteq V(X_n)$, denote by $\delta_{\text{out}} A$ the outer vertex boundary, that is the set $\{x\in V(X_n)\mid x\not\in A, x\text{ adjacent to a vertex in }A\}$.
Since $(Y_n)$ is an expander sequence, there is $h>0$ such that for all $n$ and all $B\subseteq V(Y_n)$ with $|B|\leq \frac23|Y_n|$, we have $|\delta_{\text{out}} B| \geq h|B|$.
Now let $C\subseteq V(X_n)$ with $|C|\leq \frac12|V(X_n)|$.
Write $A = C\cap V(T_n)$ and $B = C\cap V(Y_n)$.
Note that $|B|\leq \frac23|V(Y_n)|$ (provided $n$ is large enough).
We have $|\delta_{\text{out}} A| \geq \frac12|A|$, and all vertices in $\delta_{\text{out}} A$ are also in $\delta_{\text{out}} C$ unless they are in $B$, so $|\delta_{\text{out}} C|\geq \frac12|A|-|B|$.
We also have $|\delta_{\text{out}} C|\geq |\delta_{\text{out}} B\cap V(Y_n)| \geq h|B|$.
Taking a convex combination, we conclude that
\[|\delta_{\text{out}} C| \geq \left(1+\frac3{2h}\right)^{-1}\left(\frac12|A|-|B|+\frac3{2h}h|B|\right) = \frac h{2h+3}|C|.\]
Hence, $(X_n)$ is an expander sequence.

By the condition on $R_n$, we see that $(X_n)$ approximates $(Y_n)$.
By Theorem \ref{thm-approximating-prop-T}, the sequence $(X_n)$ has geometric property (T).
The $R_n$-neighbourhood of $p_n$ is the tree $T_n$, so it does not contain any cycles.
After taking a subsequence and renumbering the graphs, we get the sequence of graphs we wanted.
\end{proof}

\section*{Acknowledgement}
I would like to thank my advisor Tim de Laat for his support and suggestions.
I would also like to thank Federico Vigolo for helpful discussions.
The author is supported by the Deutsche Forschungsgemeinschaft under Germany's Excellence Strategy - EXC 2044 - 390685587, Mathematics Münster: Dynamics - Geometry - Structure.

\vspace{2cm}
\textsc{Westf\"alische Wilhelms-Universit\"at M\"unster}

\emph{Address}: Einsteinstrasse 62, 48149 M\"unster

\emph{E-mail address}: jwinkel@uni-muenster.de
\end{document}